# A Specific $M|G|\infty$ Queue System Busy Period and Busy Cycle Distributions and Parameters

**Manuel Alberto M. Ferreira**

manuel.ferreir@iscte-iul.pt

**Abstract**

Solving a Riccati equation, induced by the study of the transient behaviour of the $M|G|\infty$ queue system, a collection of service times distributions is determined. For the $M|G|\infty$ queue, which service time distribution is a member of that collection, the busy period and busy cycle probabilistic studies are performed. In extra, the properties of that distributions collection are deduced and presented.

**Keywords** Riccati Equation, Probability Distributions, $M|G|\infty$ , Busy Period, Busy Cycle

## 1. Introduction

In $M|G|\infty$ queue system the costumers arrive according to a Poisson process at rate $\lambda$. They receive a service which length is a positive random variable with distribution function G(.) and mean $\alpha$. There are always available servers. The traffic intensity is $\rho = \lambda\alpha$.

At the operation of a $M|G|\infty$ queue system, such as in any other queue system, there is a sequence of idle periods and busy periods. An idle period followed by a busy period is called a busy cycle.

When applying this queue system to real problems, the busy period and the busy cycle length probabilistic studies are very important. But their probabilistic studies are quite complex tasks. Along this work, it is shown that a more comfortable situation can be obtained solving a Riccati equation. The solution is a collection of service distributions for which both the busy period and the busy cycle have time lengths with quite simple distributions, generally related with exponential distributions and the degenerated at the origin one.

Some results for that collection of service distributions are also presented.

## 2. Riccati Equation in the *M|G|∞* Queue Occupation Study

Defining

$\beta(t) = (1-p)h(t) - \lambda p - \lambda(1-p)G(t)$, $0 \le p < 1$ ,

where $h(t) = \frac{g(t)}{1-G(t)}$ with $g(t) = G'(t)$, it is obtained

$$\frac{dG(t)}{dt} = -\lambda G^2(t) - \left(\frac{\lambda p + \beta(t)}{1-p} - \lambda\right)G(t) + \frac{\lambda p + \beta(t)}{1-p}, \quad 0 \le p < 1 \quad (1),$$

a Riccati equation in G(t), which solution is

$$G(t) = 1 - \frac{1}{\lambda}\frac{(1-e^{-\rho})e^{-\lambda t - \int_0^t \frac{\lambda p+\beta(u)}{1-p}du}}{\int_0^\infty e^{-\lambda w-\int_0^w \frac{\lambda p+\beta(u)}{1-p}du}dw - (1-e^{-\rho})\int_0^t e^{-\lambda w-\int_0^w \frac{\lambda p+\beta(u)}{1-p}du}dw},$$

$$t \ge 0, \quad -\lambda \le \frac{\int_0^t \frac{\lambda p + \beta(u)}{1-p}du}{t} \le \frac{\lambda}{e^\rho - 1},$$

$$0 \le p < 1 \quad (2),$$

see [1]. Note that, for $t \ge 0$, $G(t) = 1$ is a solution.

## 3. Busy Period and Busy Cycle Study

After (2) and considering

$$\bar{B}(s) = 1 + \lambda^{-1}\left(s - \frac{1}{\int_0^\infty e^{-st-\lambda\int_0^t[1-G(v)]dv}dt}\right) \quad (3)$$

that is the $M|G|\infty$ system busy period Laplace transform, see [2], it is obtained

$$\overline{B}(s) = \frac{1-(s+\lambda)(1-G(0))L\left[e^{-\lambda t-\int_0^t \frac{\lambda p+\beta(u)}{1-p}du}\right]}{1-\lambda(1-G(0))L\left[e^{-\lambda t-\int_0^t \frac{\lambda p+\beta(u)}{1-p}du}\right]},$$

$$-\lambda \leq \frac{\int_0^t \frac{\lambda p+\beta(u)}{1-p}du}{t} \leq \frac{\lambda}{e^\rho-1},\ 0\leq p<1 \quad (4)$$

where $L$ means Laplace transform and

$$G(0) = \frac{\lambda \int_0^\infty e^{-\lambda w-\int_0^w \frac{\lambda p+\beta(u)}{1-p}du}dw + e^{-\rho}-1}{\lambda \int_0^\infty e^{-\lambda w-\int_0^w \frac{\lambda p+\beta(u)}{1-p}du}dw}.$$

After expression (4) it is possible to compute $\frac{1}{s}\overline{B}(s)$ which inversion gives

$$B(t) = \left(1-(1-G(0))\left(\begin{array}{c}e^{-\lambda t-\int_0^t \frac{\lambda p+\beta(u)}{1-p}du}\\ +\lambda \int_0^t e^{-\lambda w-\int_0^w \frac{\lambda p+\beta(u)}{1-p}du}dw\end{array}\right)\right)*$$

$$\sum_{n=0}^{\infty} \lambda^n (1-G(0))^n \left(e^{-\lambda t-\int_0^t \frac{\lambda p+\beta(u)}{1-p}du}\right)^{*n},$$

$$-\lambda \leq \frac{\int_0^t \frac{\lambda p+\beta(u)}{1-p}du}{t} \leq \frac{\lambda}{e^\rho-1},\ 0\leq p<1 \quad (5)$$

for the $M|G|\infty$ system busy period length distribution function, being $*$ the convolution operator.

If $\beta(t)=\beta$ (constant), (2) becomes, see [3],

$$G(t) = 1 - \frac{(1-e^{-\rho})\left(\lambda+\frac{\lambda p+\beta}{1-p}\right)}{\lambda e^{-\rho}\left(e^{\left(\lambda+\frac{\lambda p+\beta}{1-p}\right)t}-1\right)+\lambda},$$

$$t\geq 0,\ -\lambda \leq \frac{\lambda p+\beta}{1-p} \leq \frac{\lambda}{e^\rho-1},\ 0\leq p<1 \quad (6)$$

and (5),

$$B^\beta(t) = 1 - \frac{\lambda+\frac{\lambda p+\beta}{1-p}}{\lambda}(1-e^{-\rho})\ e^{-e^{-\rho}\left(\lambda+\frac{\lambda p+\beta}{1-p}\right)t},$$

$$t\geq 0, -\lambda \leq \frac{\lambda p+\beta}{1-p} \leq \frac{\lambda}{e^\rho-1}, 0\leq p<1 \quad (7).$$

So, if the service distribution is given by (6), the busy period distribution function is the mixture of a degenerated distribution at the origin with an exponential.

For $\beta=-\lambda,\ 0\leq p<1$,

$$G(t)=1,\ t\geq 0 \quad \text{and}\ B^\beta(t)=1,\ t\geq 0$$

that is if the service distribution is degenerated at the origin with probability 1 it happens the same with the busy period distribution.

Finally note that for $\beta=\lambda\frac{1-pe^\rho}{e^\rho-1},\ 0\leq p<1$

$$B^\beta(t) = 1-e^{-\frac{\lambda}{e^\rho-1}t},\ t\geq 0$$

(Purely exponential). And $B(t)$ given by (5) satisfies

$$B(t)\geq 1-e^{-\frac{\lambda}{e^\rho-1}t},\ t\geq 0.$$

Calling $\overline{Z}(s)$ the busy cycle length Laplace transform, and noting that the $M|G|\infty$ queue idle period is exponentially distributed with parameter $\lambda$ - as it happens for any queue system with Poisson arrivals at rate $\lambda$ - and that the idle period and the busy period are independent for the $M|G|\infty$ queue, see [4], it is obtained

$$\overline{Z}(s) = \frac{\lambda}{\lambda+s}\overline{B}(s) \quad (8).$$

After (4) and (8) it is possible to compute $(1/s)\overline{Z}(s)$ which inversion gives

$$Z(t) = (\lambda e^{-\lambda t})$$
$$*\left(\begin{array}{c}1-(1-G(0))\\ \left(\left(e^{-\lambda t-\int_0^t \frac{\lambda p+\beta}{1-p}du}+\lambda\int_0^t e^{-\lambda w-\int_0^w \frac{\lambda p+\beta}{1-p}du}dw\right)\right)\end{array}\right)*$$

$$*\sum_{n=0}^{\infty}\lambda^n(1-G(0))^n\left(e^{-\lambda t-\int_0^t \frac{\lambda p+\beta}{1-p}du}\right)^{*n},$$

$$-\lambda \leq \frac{\int_0^t \frac{\lambda p+\beta}{1-p}du}{t} \leq \frac{\lambda}{e^\rho-1}, 0\leq p<1 \quad (9)$$

for the busy cycle distribution function.

If $\beta(t)=\beta$ (constant) (9) becomes

$$Z^\beta(t) = 1 - \frac{(1-e^{-\rho})\left(\lambda+\frac{\lambda p+\beta}{1-p}\right)}{\lambda-e^{-\rho}\left(\lambda+\frac{\lambda p+\beta}{1-p}\right)}\ e^{-e^{-\rho}\left(\lambda+\frac{\lambda p+\beta}{1-p}\right)t}+$$

$$\frac{\frac{\lambda p+\beta}{1-p}}{\lambda-e^{-\rho}\left(\lambda+\frac{\lambda p+\beta}{1-p}\right)}\ e^{-\lambda t},$$

$$t\geq 0,\ -\lambda \leq \frac{\lambda p+\beta}{1-p} \leq \frac{\lambda}{e^\rho-1},\ 0\leq p<1 \quad (10).$$

So, if the service time distribution function is given by (6) the busy cycle distribution function is the mixture of two exponential distributions.

Note that, for $\beta=\lambda\frac{1-pe^\rho}{e^\rho-1},\ 0\leq p<1$

$$Z^\beta(t) = 1 - \frac{(e^\rho-1)\ e^{-\frac{\lambda}{e^\rho-1}t}-e^{-\lambda t}}{e^\rho-2},\ t\geq 0.$$

And $Z(t)$, given by (9), satisfies

$$Z(t) \geq 1 - \frac{(e^\rho-1)e^{-\frac{\lambda}{e^\rho-1}t} - e^{-\lambda t}}{e^\rho - 2}, \quad t \geq 0 \quad (11).$$

Note that (11) is coherent even for $\rho = \ln 2$ since

$$\lim_{\rho \to \log 2} \left(1 - \frac{(e^\rho-1)e^{-\frac{\lambda}{e^\rho-1}t} - e^{-\lambda t}}{e^\rho - 2}\right) =$$

$$\lim_{\rho \to \log 2} \frac{e^\rho - 2 - (e^\rho-1)e^{-\frac{\lambda}{e^\rho-1}t} - e^{-\lambda t}}{e^\rho - 2} = 1 - (1+\lambda t)e^{-\lambda t}.$$

For $\beta = -\lambda$, $0 \leq p < 1$

$$Z^\beta(t) = 1 - e^{-\lambda t}, \ t \geq 0.$$

Finally, also for $Z(t)$ given by (9),

$$Z(t) \leq 1 - e^{-\lambda t}, \quad t \geq 0 \quad (12).$$

## 4. Determining Busy Period and Busy Cycle Parameters

As for the busy period moments, calling $B$ the busy period random variable, note that (3) is equivalent, see [5], to

$$(\overline{B}(s) - 1)(C(s) - 1) = \lambda^{-1} s C(s),$$

being

$$C(s) = \int_0^\infty e^{-st - \lambda \int_0^t [1-G(v)]dv} \lambda(1 - G(t))dt.$$

Differentiating $n$ times, using Leibnitz's formula, and making $s = 0$, it results

$$E[B^n] = (1)^{n+1} \left\{ \frac{e^\rho}{\lambda} n C^{(n-1)}(0) - e^\rho \sum_{p=1}^{n-1} (-1)^{n-p} \binom{n}{p} \cdot E[B^{n-p}] C^{(p)}(0) \right\},$$
$$n = 1, 2, \ldots \quad (13)$$

where

$$C^{(n)}(0) = \int_0^\infty (-t)^n e^{-\lambda \int_0^t [1-G(v)]dv} \lambda(1 - G(t))dt,$$
$$n = 0,1,2, \quad (14).$$

The expression (13) gives a recursive method to compute $E[B^n], n = 1,2,\ldots$ as a function of $C^{(n)}(0), n = 0,1,2,\ldots$

So, it is possible to write exact expressions for the moments centred at the origin.

Making $n = 1$ it is obtained

$$E[B] = \frac{e^\rho - 1}{\lambda} \quad (15)$$

as it is well known, independently of the service time distribution.

For the busy cycle moments, if Z is the busy cycle random variable,

$$E[Z^n] = \sum_{p=0}^{n} \binom{n}{p} E[B^p] \cdot \frac{(n-p)!}{\lambda^{n-p}} \quad (16).$$

The $M|G|\infty$ queue busy period "peak" is the Laplace transform value at $1/\alpha$, see [6]. It is a parameter that characterizes the busy period distribution and contains information about all its moments.

For the collection of service distributions (2) the "peak", called $pi$, is

$$pi = \frac{1 - \left(\frac{1}{\alpha} + \lambda\right)(1 - G(0)) L\left[e^{-\lambda t - \int_0^t \frac{\lambda p + \beta(u)}{1-p} du}\right]\left(\frac{1}{\alpha}\right)}{1 - \lambda(1 - G(0)) L\left[e^{-\lambda t - \int_0^t \frac{\lambda p + \beta(u)}{1-p} du}\right]\left(\frac{1}{\alpha}\right)},$$

$$-\lambda \leq \frac{\int_0^t \frac{\lambda p + \beta(u)}{1-p} du}{t} \leq \frac{\lambda}{e^\rho - 1}, \ 0 \leq p < 1 \quad (17).$$

For $\beta(t) = \beta$ (constant), see [7],

$$pi = \frac{e^{-\rho}(\lambda + \beta)(\rho + 1) - \lambda p - \beta}{\lambda(e^{-\rho}(\rho + \alpha\beta) + 1 - p)},$$
$$-\lambda \leq \beta \leq \frac{\lambda(1-pe^\rho)}{e^\rho - 1}, 0 \leq p < 1 \quad (18)$$

In [6] it is also introduced another measure: the "modified peak" got after the "peak" taking out the terms that are permanent for the busy period in different service distributions and putting over the common part. Call it $qi$. It is easy to verify that $qi = pi \frac{\rho}{e^\rho - \rho - 1} + 1$ and, so, for the distributions given by collection (2) with $\beta(t) = \beta$ (constant), so collection (6), see again [7],

$$qi = \frac{e^{-\rho}(\lambda + \beta)(\rho + 1) - \lambda p + \beta}{\lambda(e^{-\rho}(\rho + \alpha\beta) + 1 - p)} \cdot \frac{\rho}{e^\rho - \rho - 1} + 1,$$
$$-\lambda \leq \beta \leq \frac{\lambda(1-pe^\rho)}{e^\rho - 1}, 0 \leq p < 1 \quad (19)$$

For the busy cycle of the $M|G|\infty$ queue the "peak", see [8], called now $pi'$, for the service distributions given by the collection (2) is

$$pi' = \frac{\rho}{\rho+1} \frac{1-\left(\frac{1}{\alpha}+\lambda\right)(1-G(0))L\left[e^{-\lambda t-\int_0^t \frac{\lambda p+\beta(u)}{1-p}du}\right]\left(\frac{1}{\alpha}\right)}{1-\lambda(1-G(0))L\left[e^{-\lambda t-\int_0^t \frac{\lambda p+\beta(u)}{1-p}du}\right]\left(\frac{1}{\alpha}\right)},$$

$$-\lambda \leq \frac{\int_0^t \frac{\lambda p+\beta(u)}{1-p}du}{t} \leq \frac{\lambda}{e^\rho-1}, \quad 0 \leq p < 1 \quad (20).$$

With $\beta(t) = \beta$ (constant), see [7],

$$pi' = \alpha \frac{e^{-\rho}(\lambda+\beta)(\rho+1)-\lambda p-\beta}{(\rho+1)(e^{-\rho}(\rho+\alpha\beta)+1-p)},$$
$$-\lambda \leq \beta \leq \frac{\lambda(1-pe^\rho)}{e^\rho-1}, 0 \leq p < 1 \quad (21).$$

And the "modified peak", called now $qi'$, $pi'\frac{\rho}{e^\rho-\rho}+1$, for the service distributions given by the collection (6) is

$$qi' = \alpha \frac{e^{-\rho}(\lambda+\beta)(\rho+1)-\lambda p-\beta}{(\rho+1)(e^{-\rho}(\rho+\alpha\beta)+1-p)} \frac{\rho}{e^\rho-\rho}+1,$$
$$-\lambda \leq \beta \leq \frac{\lambda(1-pe^\rho)}{e^\rho-1}, 0 \leq p < 1 \quad (22).$$

The busy cycle renewal function value of the $M|G|\infty$ queue, at the instant $t$, gives the mean number of busy periods that begin in $[0, t]$, see [9]. Calling it $R(t)$:

$$R(t) = e^{-\lambda \int_0^t [1-G(v)]dv} + \lambda \int_0^t e^{-\lambda \int_0^u [1-G(v)]dv} du \quad (23).$$

If the service time is a random variable with distribution function given by a member of the collection (6)

$$R(t) = e^{-\rho}(1+\lambda t) + (1-e^{-\rho})\frac{\lambda p+\beta}{\lambda+\beta}e^{-\left(\lambda+\frac{\lambda p+\beta}{1-p}\right)t} +$$
$$(1-e^{-\rho})\frac{\lambda p+\beta}{\lambda+\beta}, -\lambda \leq \beta \leq \frac{\lambda(1-pe^\rho)}{e^\rho-1}, 0 \leq p < 1 \quad (24).$$

## 5. Collection (2) Distributions Moments Determination

If, in (2), $G_i(t)$ is the solution associated to $\rho_i$, $i = 1,2,3,4$ it is easy to see that

$$\frac{G_4(t)-G_2(t)}{G_4(t)-G_1(t)} \cdot \frac{G_3(t)-G_1(t)}{G_3(t)-G_2(t)} = \frac{e^{-\rho_4}-e^{-\rho_2}}{e^{-\rho_4}-e^{-\rho_1}} \cdot \frac{e^{-\rho_3}-e^{-\rho_1}}{e^{-\rho_3}-e^{-\rho_2}} \quad (25)$$

as it had to happen since the differential equation considered is a Riccati one.

And computing,

$$\int_0^\infty [1-G(t)]\, dt =$$

$$\int_0^\infty \frac{1}{\lambda} \frac{(1-e^{-\rho})C(t)}{\int_0^\infty e^{-\lambda w-\int_0^w \frac{\lambda p+\beta(u)}{1-p}du}dw - (1-e^{-\rho})\int_0^t C(w)dw} dt$$
$$= \alpha,$$

where $C(v) = e^{-\lambda v-\int_0^v \frac{\lambda p+\beta(u)}{1-p}du}, v = t, w$ in accordance with the fact that (2) is a collection of positive random variables distribution functions.

The density associated to $G(t)$ given by (2) is

$$g(t) = -\frac{1}{\lambda}\left[\frac{A''(t)}{B-A(t)} + \frac{[A^t(t)]^2}{[B-A(t)]^2}\right],$$
$$t \geq 0, \quad -\lambda \leq \frac{\int_0^t \frac{\lambda p+\beta(u)}{1-p}du}{t} \leq \frac{\lambda}{e^\rho-1},$$
$$0 \leq p < 1 \quad (26)$$

where

$$A(t) = (1-e^{-\rho})\int_0^t e^{-\lambda w-\int_0^w \frac{\lambda p+\beta(u)}{1-p}du}dw$$

and

$$B = \int_0^\infty e^{-\lambda w-\int_0^w \frac{\lambda p+\beta(u)}{1-p}du}dw.$$

If $\beta(t) = \beta$ (constant), see [7],

$$g(t) = \frac{(1-e^{-\rho})e^{-\rho}\left(\lambda+\frac{\lambda p+\beta}{1-p}\right)^2 e^{-\left(\lambda+\frac{\lambda p+\beta}{1-p}\right)t}}{\lambda\left[e^{-\rho}+(1-e^{-\rho})e^{-\left(\lambda+\frac{\lambda p+\beta}{1-p}\right)t}\right]^2},$$

$$t > 0, -\lambda \leq \beta \leq \frac{\lambda(1-pe^{-\rho})}{e^\rho-1}, 0 \leq p < 1 \quad (27)$$

So,

$$\int_0^\infty t^n g(t)dt = \frac{(1-e^{-\rho})e^{-\rho}\left(\lambda+\frac{\lambda p+\beta}{1-p}\right)^2}{\lambda} \times$$
$$\int_0^\infty t^n \frac{e^{-\left(\lambda+\frac{\lambda p+\beta}{1-p}\right)t}}{\left[e^{-\rho}+(1-e^{-\rho})e^{-\left(\lambda+\frac{\lambda p+\beta}{1-p}\right)t}\right]^2} dt \quad (28)$$

But, $\int_0^\infty t^n \frac{e^{-\left(\lambda+\frac{\lambda p+\beta}{1-p}\right)t}}{\left[e^{-\rho}+(1-e^{-\rho})e^{-\left(\lambda+\frac{\lambda p+\beta}{1-p}\right)t}\right]^2} dt$

$$\geq \int_0^\infty t^n e^{-\left(\lambda+\frac{\lambda p+\beta}{1-p}\right)t} dt =$$

$$\frac{1}{\lambda+\frac{\lambda p+\beta}{1-p}} \frac{n!}{\left(\lambda+\frac{\lambda p+\beta}{1-p}\right)^n}, \beta \neq -\lambda.$$

$$\int_0^\infty t^n \frac{e^{-\left(\lambda + \frac{\lambda p + \beta}{1-p}\right)t}}{\left[e^{-\rho} + (1-e^{-\rho})e^{-\left(\lambda + \frac{\lambda p + \beta}{1-p}\right)t}\right]^2} dt \leq$$

$$e^{2\rho} \int_0^\infty t^n e^{-\left(\lambda + \frac{\lambda p + \beta}{1-p}\right)t} dt =$$

$$\frac{e^{2\rho}}{\lambda + \frac{\lambda p + \beta}{1-p}} \frac{n!}{\left(\lambda + \frac{\lambda p + \beta}{1-p}\right)^n}, \beta \neq -\lambda.$$

So, calling $T$ the random variable corresponding to $G(t)$:

$$\frac{(1-e^{-\rho})e^{-\rho}}{\lambda} \frac{n!}{\left(\lambda + \frac{\lambda p + \beta}{1-p}\right)^{n-1}} \leq E[T^n]$$

$$\leq \frac{e^\rho - 1}{\lambda} \frac{n!}{\left(\lambda + \frac{\lambda p + \beta}{1-p}\right)^{n-1}},$$

$$-\lambda < \beta \leq \frac{\lambda(1-pe^{-\rho})}{e^\rho - 1}, 0 \leq p < 1, n = 1,2,. \quad (29)$$

**Notes:**

-The expression (29), giving bounds for $E[T^n]$, guarantees its existence,

-For $n = 1$ the expression (29) is useless because $E[T] = \alpha$. Note, curiously, that the upper bound is $\frac{e^\rho - 1}{\lambda}$, the $M/G/\infty$ system busy period mean value,

-For $n = 2$, subtracting to both bounds, $\alpha^2$ from expression (29) result bounds for $VAR[T]$,

-For
$$\beta = -\lambda, \ E[T^n] = 0, n = 1,2,\ldots,$$

evidently.

See, however, that (6) can be put like:

$$G(t) = \frac{1 + \frac{\lambda p + \beta}{1-p}(1-e^\rho)e^{-\left(\lambda + \frac{\lambda p + \beta}{1-p}\right)t}}{1 - (1-e^\rho)e^{-\left(\lambda + \frac{\lambda p + \beta}{1-p}\right)t}},$$

$$t \geq 0, -\lambda \leq \beta \leq \frac{\lambda(1-pe^\rho)}{e^\rho - 1}, 0 \leq p < 1 \quad (30)$$

And, for $\rho < \ln 2$,

$$G(t) = \left(1 + \frac{\frac{\lambda p + \beta}{1-p}}{\lambda}(1-e^\rho)e^{-\left(\lambda + \frac{\lambda p + \beta}{1-p}\right)t}\right) \cdot$$

$$\sum_{k=0}^\infty (1-e^\rho)^k e^{-k\left(\lambda + \frac{\lambda p + \beta}{1-p}\right)t},$$

$$t \geq 0, -\lambda \leq \beta \leq \frac{\lambda(1-pe^\rho)}{e^\rho - 1}, 0 \leq p < 1 \quad (31)$$

After (31) the $T$ Laplace Transform for $\rho < \ln 2$, can be easily derived. And, so

- For $\rho < \ln 2$

$$E[T^n] = -\left(1 + \frac{\frac{\lambda p + \beta}{1-p}}{\lambda}\right) n! \sum_{k=1}^\infty \frac{(1-e^\rho)^k}{k\left(\lambda + \frac{\lambda p + \beta}{1-p}\right)^n},$$

$$-\lambda < \beta \leq \frac{\lambda(1-pe^\rho)}{e^\rho - 1}, 0 \leq p < 1, n = 1,2,\ldots \quad (32)$$

**Notes:**

-For $n = 1$,

$$E[T] = -\left(1 + \frac{\frac{\lambda p + \beta}{1-p}}{\lambda}\right) \sum_{k=1}^\infty \frac{(1-e^\rho)^k}{k\left(\lambda + \frac{\lambda p + \beta}{1-p}\right)} =$$

$$\frac{1}{\lambda} \sum_{k=1}^\infty (-1)^{k+1} \frac{(1-e^\rho)}{k} = \frac{1}{\lambda} \log e^\rho = \frac{\rho}{\lambda} = \alpha$$

-For $n \geq 2$ only a finite number of parcels in the infinite sum must be taken. Calling $M$ this number, to get an error lesser than $\varepsilon$ it is necessary to have simultaneously

$$M > \frac{1}{\lambda + \frac{\lambda p + \beta}{1-p}} - 1,$$

$$M > \log_{(e^\rho - 1)} \frac{\varepsilon e^\rho \lambda}{n!\left(\lambda + \frac{\lambda p + \beta}{1-p}\right)} - 1.$$

So, it is evident now that this distributions collection moments computation is a complex task. This was already true for the study of [10] where the results presented are a particular situation of these ones for $p = 0$.

The consideration of the approximation

$$E_m^n = \sum_{k=1}^\infty \left(\frac{k}{m}\right)^n \left[G\left(\frac{k}{m}\right) - G\left(\frac{k-1}{m}\right)\right],$$

$$-\lambda < \beta \leq \frac{\lambda(1-pe^\rho)}{e^\rho - 1}, 0 \leq p < 1, n = 1,2,\ldots \quad (33)$$

may be helpful since $\lim_{m \to \infty} E_m^n = E[T^n], n = 1,2,\ldots$, see [11], that allow the numerical computation of the moments.

# 6. Conclusions

The function $\beta(t)$, defined in 1., that leads to the equation (1), the source of the results presented, is induced by the $M/G/\infty$ transient probabilities monotony study, see [1] and [3]. Being well known the close relations existing among the transient behaviour and the busy period distributions, it is not surprising that the service distribution functions, solutions of (1), give rise to so friendly distributions to the $M/G/\infty$ queue busy period. This is a quite unusual situation, that allows simple computations for the busy period and busy cycle distributions probabilities and parameters.

This so simple structure lies on the exponential distribution and the deterministic one. Among the cases seen there are situations of purely deterministic and exponential distributions, mixtures of deterministic and exponential distributions and mixtures of two exponential distributions. The great presence of the exponential distribution in this context is due certainly to the so friendly properties of the Poisson Process and to the structure of the function resulting from the inversion of the Laplace Transform given by (4).

Even the stochastic processes related with this queue become quite simple. They are very close to the Poisson Process and in certain situations are even Poisson Processes.

In conclusion, all of this results in a very simple expression for the inversion of (4), with the service distributions solutions of (1). In consequence, as it was shown in this work, the distribution functions, and the parameters for the busy period and, consequently, for the busy cycle result very simple to compute and interpret.

Many of these results are true for other service distribution functions, in the case of insensibility: parameters that depend on the service distribution only through its mean. Some others may result in good approximations.

Note still that when the distribution functions for $B(t)$ and $Z(t)$ are a little more complicated, simple bounds are presented.

In the field of Queuing Theory this situation seems unique, there having no news of anything similar.

For more details on this subject see [12, 13, and 14].

# Acknowledgements

This work is partially financed by national funds through FCT - Fundação para a Ciência e Tecnologia, I.P., under the project FCT UIDB/04466/2020. Furthermore, the author thanks the ISCTE-IUL and ISTAR-IUL, for their support.